\documentclass[11pt]{amsart}
\usepackage{mathrsfs}
\usepackage[active]{srcltx}
\usepackage{mathrsfs}

\newtheorem{proposition}{Proposition}[section]
\newtheorem{theorem}{Theorem}[section]

\newtheorem{corollary}{Corollary}[section]

\newtheorem{remark}{Remark}[section]

\numberwithin{equation}{section}

\begin{document}

\title[Sharp Morrey-Sobolev inequalities on Riemannian Manifolds]{Sharp Morrey-Sobolev inequalities on complete Riemannian Manifolds}

\author{Alexandru Krist\'aly}
\address{Department of Economics, Babe\c s-Bolyai University, 400591 Cluj-Napoca,
Romania,} \address{Institute of Applied Mathematics, \'Obuda
University, 1034 Budapest, Hungary\\ }
\email{alexandrukristaly@yahoo.com}

\thanks{Research supported by a grant of the Romanian National Authority for Scientific Research,
CNCS-UEFISCDI,  "Symmetries in elliptic problems: Euclidean and
non-Euclidean techniques", project no. PN-II-ID-PCE-2011-3-0241, and
by the J\'anos Bolyai Research Scholarship of the Hungarian Academy
of Sciences. The present work was initiated during the author's
visit at the Institut des Hautes \'Etudes Scientifiques (IH\'ES),
Bures-sur-Yvette, France.}

\subjclass[2000]{Primary 58J60; Secondary 53C21}



\keywords{Morrey-Sobolev inequality; Riemannian manifold; sharp
constant; curvature; rigidity}

\begin{abstract}
Two Morrey-Sobolev inequalities (with support-bound and $L^1-$bound,
respectively) are investigated on complete Riemannian manifolds with
their sharp constants in $\mathbb R^n$. We prove the following
results in both cases:

$\bullet$ If $(M,g)$ is a {\it Cartan-Hadamard manifold} which
verifies the $n-$di\-men\-sional Cartan-Hadamard conjecture,  sharp
Morrey-Sobolev inequalities hold on $(M,g)$. Moreover, extremals
exist if and only if $(M,g)$ is isometric to the standard Euclidean
space $(\mathbb R^n,e)$.

$\bullet$ If $(M,g)$ has {\it non-negative Ricci curvature}, $(M,g)$
supports the sharp Morrey-Sobolev inequalities if and only if
$(M,g)$ is isometric to $(\mathbb R^n,e)$.

\end{abstract}

\maketitle

 \section{Introduction and main results}

One of the most important topics of Sobolev inequalities is to find
sharp constants and extremals in the embeddings $W^{1,p}(\mathbb
R^n)\hookrightarrow L^q(\mathbb R^n)$, where the numbers $p,q\in
\mathbb R$ and $n\in \mathbb N$ are related in the Sobolev sense.
Owing to a systematic study initiated by T. Aubin and G. Talenti in
the middle of seventies, various results are available nowadays
concerning sharp constants and extremals in Sobolev inequalities
both in the Euclidean and Riemannian frameworks; see  Ghoussoub and
Moradifam \cite{Gh-Mor-monograph}, Hebey \cite{Hebey}, Maz'ya
\cite{Mazja}, and references therein. We emphasize that sharp
Sobolev inequalities in $\mathbb R^n$ were mostly studied for $ p
\in [1, n)$; see the pioneering works of Federer and Fleming
\cite{FF} when $p = 1,$ and Aubin \cite{Aubin} and Talenti
\cite{Talenti-masik} when $1 < p < n.$ Moreover, when $ p \in [1,
n)$, several rigidity results can be found on Riemannian manifolds
supporting Sobolev-type inequalities with their Euclidean sharp
constants, see Ledoux \cite{Ledoux}, do Carmo and Xia
\cite{doCarmo-Xia}, Druet and Hebey \cite{DH}, Druet, Hebey and
Vaugon \cite{Druetetal}, and the comprehensive monograph of Hebey
\cite{Hebey}.

The main purpose of this paper is to investigate  {\it sharp
Morrey-Sobolev inequalities} (i.e., $p>n$) on non-compact complete
Riemannian manifolds having either non-positive sectional curvature
or non-negative Ricci curvature.  Hereafter, sharpness means that a
given inequality on the Riemannian manifold is valid with its
Euclidean sharp constant.

The Morrey-Sobolev inequality in $\mathbb R^n$ states that the
embedding $W^{1,p}(\mathbb R^n)\hookrightarrow L^\infty(\mathbb
R^n)$ is continuous when $p>n$ (see \cite{Brezis}, \cite{Mazja}),
i.e., there exists ${\sf C}(p,n)>0$ such that
$$\|u\|_{L^\infty(\mathbb R^n)}\leq {\sf C}({p,n})(\|u\|_{L^p(\mathbb
R^n)}+\|\nabla u\|_{L^p(\mathbb R^n)}),\ \forall u\in
W^{1,p}(\mathbb R^n).$$  A generic Morrey-Sobolev inequality has
been established on smooth complete $n-$di\-men\-sional Riemannian
manifolds with Ricci curvature bounded from below and verifying a
global volume growth assumption, see Coulhon \cite{Coulhon}.

Let $(M,g)$ be an $n(\geq 2)-$di\-men\-sional smooth complete
Riemannian manifold.  In order to present our results, we need two
notions. First, we say that a function $u:M\to [0,\infty)$ is {\it
concentrated around} $x_0\in M$ if for every
$0<t<\|u\|_{L^\infty(M)},$ the level set $\{x\in M:u(x)>t\}$ is a
geodesic ball $B(x_0,\rho_t)=\{x\in M:d(x_0,x)<\rho_t\}$ for some
$\rho_t>0$. Hereafter, $d:M\times M\to \mathbb R$ denotes the usual
distance function associated with $g$. Second, Morrey-Sobolev
inequalities will be particularly investigated on Car\-tan-Hadamard
manifolds as well (i.e., on simply connected, complete Riemannian
manifolds with non-positive sectional
curvature) where the validity of the {\it Cartan-Hadamard conjecture} will play an indispensable role. For the sake of completeness,  we recall the\\

\noindent {\bf Cartan-Hadamard conjecture in $n-$dimension} (see
Aubin \cite{Aubin}). {\it Let $(M,g)$ be an $n-$dimen\-sional
Cartan-Ha\-da\-mard manifold. Then any compact domain $D\subset M$
with smooth boundary $\partial D$ satisfies the Euclidean
isoperimetric inequality, i.e.,
\begin{equation}\label{cartan-hadamard-conjecture}
    {\rm Area}_g(\partial D) \geq  n\omega_n^\frac{1}{n} {\rm
    Vol}_g(D)^\frac{n-1}{n}.
\end{equation}
Moreover, equality holds in {\rm (\ref{cartan-hadamard-conjecture})}
if and only if $D$ is isometric to the $n-$dimensional Euclidean
ball with volume ${\rm Vol}_g(D)$.}\\

\noindent Hereafter, $\omega_n$ is the volume of the
$n-$di\-men\-sional Euclidean unit ball; Area$_g(\partial D)$ stands
for the area of $\partial D$ with respect to the metric induced on
$\partial D$ by $g$; and Vol$_g(D)$ is the volume of $D$ with
respect to $g$.

\begin{remark}\rm Cartan-Hadamard conjecture is true
in dimension $2,$ see Weil \cite{Weil}; in dimension $3,$ see
Kleiner \cite{Kleiner};  and in dimension $4$, see Croke
\cite{Croke}, but it is open for higher dimensions.
\end{remark}

Now we are ready to present our main results.\\

\noindent (I) {\bf Sharp Morrey-Sobolev inequality with
support-bound.}  
Let $(M,g)$ be an $n(\geq 2)-$di\-men\-sional smooth complete
Riemannian manifold and $p>n$.  For some $C>0,$ we consider on
$(M,g)$ the Morrey-Sobolev inequality
 {$$
    \| u\|_{L^\infty(M)} \leq C \mathcal H^n({\rm sprt}\ u )^{\frac{1}{n}-\frac{1}{p}} \| \nabla_gu\|_{L^p(M)},\ \ \forall u\in {\rm Lip}_0(M). \eqno{{\bf ({MS})}_{C}^1}
$$}
Here, ${\rm sprt}\ u$ is the support of $u$, $\mathcal H^n$ is the
$n-$dimensional Hausdorff  measure on $M$, $\|\nabla_g u\|_{L^p(M)}$
stands for the $L^p(M)$ norm of the vector  $\nabla_g u(x)\in T_xM$,
while ${\rm Lip}_0(M)$ is the space of Lipschitz functions with
compact support defined on $M$. Although we can put $C_0^\infty(M)$
instead of
 ${\rm Lip}_0(M)$ in ${\bf ({MS})}_{C}^1$ due to density reasons, we prefer the latter
 choice taking into account the specific shape of extremals
 in the Euclidean setting.
Indeed, by using symmetrization and rearrangement arguments, Talenti
\cite[Theorem 2.E]{Talenti} proved that if $(M,g)=(\mathbb R^n,e)$
is the standard Euclidean space, then ${\bf
({MS})}_{{\textsf{C}_1(p,n)}}^1$ holds with the sharp constant
\begin{equation}\label{C1-konstans}
    \textsf{C}_1(p,n)=n^{-\frac{1}{p}}\omega_n^{-\frac{1}{n}}\left(\frac{p-1}{p-n}\right)^\frac{1}{p'},
\end{equation}
where $p'=\frac{p}{p-1},$ and the unique class of extremals (up to a
constant multiplication) is given by
\begin{equation}\label{extremal-1}
    u_{\lambda,x_0}(x)=\left(\lambda^\frac{p-n}{p-1}-|x-x_0|^\frac{p-n}{p-1}\right)_+,\
x\in \mathbb R^n,
\end{equation}
where $\lambda>0$, $x_0\in \mathbb R^n,$ and $r_+=\max(r,0).$
 Clearly, for a fixed $x_0\in \mathbb R^n$, the function
$u_{\lambda,x_0}\in {\rm Lip}_0(\mathbb R^n)$ is concentrated around
$x_0$, its support being the closed Euclidean ball
$B_e[x_0,\lambda]$.

Our first result reads as follows.

\begin{theorem}\label{theorem-morrey-1} Let $(M,g)$ be
an $n-$di\-men\-sional Car\-tan-Hadamard manifold which verifies the
Cartan-Hadamard conjecture in the same dimension, and let $p>n.$
\begin{itemize}
  \item[{\rm (i)}] {\rm [Sharpness]} The Morrey-Sobolev inequality ${\bf
({MS})}_{{\sf C}_1(p,n)}^1$ holds on $(M,g);$ moreover, ${\sf
C}_1(p,n)$ is sharp, i.e.,{
$${\sf C}_1(p,n)^{-1}=\inf_{u\in {\rm Lip}_0(M)\setminus \{0\}} \frac{\mathcal H^n({\rm sprt}\ u )^{\frac{1}{n}-\frac{1}{p}} \| \nabla_gu\|_{L^p(M)}}{\|
u\|_{L^\infty(M)}}.$$}
  \item[{\rm (ii)}] {\rm [Extremals]} Let $x_0\in M$. The following
statements are equivalent:
\begin{itemize}
  \item[\rm (a)] For every $\kappa>0$ there exists a non-negative extremal function $u\in {\rm Lip}_0(M)$
  in ${\bf ({MS})}_{{\sf C}_1(p,n)}^1$, concentrated around $x_0$
  and $\mathcal H^n({\rm sprt}\ u )=\kappa;$
   \item[\rm (b)] $(M,g)$ is isometric to $(\mathbb R^n,e)$.
\end{itemize}
\end{itemize}
\end{theorem}

\noindent In the  non-negatively curved case we state the following
result:

\begin{theorem}\label{theorem-morrey-2} Let $(M,g)$ be a  complete,
 $n-$di\-men\-sional Riemannian manifold with non-negative Ricci
 curvature, let $p>n,$ and assume that ${\bf ({MS})}_{C}^1$ holds on $(M,g)$
 for some $C>0.$ Then the following assertions hold:
 \begin{itemize}
   \item[\rm (i)] $C\geq {\sf C}_1(p,n);$
   \item[\rm (ii)] $(M,g)$ has the large volume balls property, i.e., $${\rm
   Vol}_g(B(x,\rho))\geq \left(\frac{{\sf C}_1(p,n)}{C}\right)^\frac{pn}{p-n}\omega_n\rho^n\ {for\ all}\ x\in M,\rho\geq 0;$$
   \item[\rm (iii)] ${\bf ({MS})}_{{\sf C}_1(p,n)}^1$ holds on
   $(M,g)$ if and only if $(M,g)$ is isometric to $(\mathbb R^n,e).$
 \end{itemize}
\end{theorem}

\vspace{0.3cm}

\noindent (II) {\bf Sharp Morrey-Sobolev inequality with
$L^1$-bound.} Instead of having a support-bound estimate in term of
$\mathcal H^n({\rm sprt}\ u )$ for $\|u\|_{L^\infty(M)}$, we can use
a suitable interpolation between $\| \nabla_gu\|_{L^p(M)}$ and some
other norm $\|u\|_{L^q(M)}$, $q\in
[1,\infty)$. To do this, let $(M,g)$ be a smooth $n-$dimensional complete Riemannian manifold and $p>n$. 
For some $C>0,$ we consider on $(M,g)$ the Morrey-Sobolev inequality
 {$$
    \| u\|_{L^\infty(M)} \leq C \|u\|_{L^1(M)}^{1-\eta} \| \nabla_gu\|_{L^p(M)}^\eta,\ \ \forall u\in {\rm Lip}_0(M), \eqno{{\bf ({MS})}_{C}^2}
$$}
where
\begin{equation}\label{eta-def}
    \eta=\frac{np}{np+p-n}.
\end{equation}
 Talenti \cite[Theorem
2.C]{Talenti} proved that ${\bf ({MS})}_{{\sf C}_2(p,n)}^2$ holds on
$(\mathbb R^n,e)$ with the sharp constant
$${\sf
C}_2(p,n)=(n\omega_n^\frac{1}{n})^{-\frac{np'}{n+p'}}\left(\frac{1}{n}+\frac{1}{p'}\right)
\left(\frac{1}{n}-\frac{1}{p}\right)^\frac{(n-1)p'-n}{n+p'}
\left({\sf B}\left(
\frac{1-n}{n}p'+1,p'+1\right)\right)^\frac{n}{n+p'},$$ where ${\sf
B}(\cdot,\cdot)$ stands for the Euler beta-function. The unique
family of extremals (up to a constant multiplication) is given by
$$v_{\lambda,x_0}(x)=\left\{
\begin{array}{lll}
\displaystyle\int_{|x-x_0|}^\lambda
r^\frac{1-n}{p-1}(\lambda^n-r^n)^\frac{1}{p-1}{\rm d}r, & & {\rm
if}\ |x-x_0|\leq \lambda;
\\ 0,&  & {\rm otherwise},
\end{array}
\right.$$ where $\lambda>0$, $x_0\in \mathbb R^n.$

Similar results can be obtained for ${\bf ({MS})}_{{\sf
C}_2(p,n)}^2$ as in Theorems \ref{theorem-morrey-1} \&
\ref{theorem-morrey-2}; namely, we prove:

\begin{theorem}\label{theorem-morrey-3}
Let $(M,g)$ be an $n-$di\-men\-sional Car\-tan-Hadamard manifold
which verifies the Cartan-Hadamard conjecture in the same dimension,
and let $p>n.$
\begin{itemize}
  \item[{\rm (i)}] {\rm [Sharpness]} The Morrey-Sobolev inequality ${\bf
({MS})}_{{\sf C}_2(p,n)}^2$ holds on $(M,g);$ moreover, ${\sf
C}_2(p,n)$ is sharp, i.e.,{
$${\sf C}_2(p,n)^{-1}=\inf_{u\in {\rm Lip}_0(M)\setminus \{0\}} \frac{\|u\|_{L^1(M)}^{1-\eta} \| \nabla_gu\|_{L^p(M)}^\eta}{\|
u\|_{L^\infty(M)}},$$} where $\eta$ is given by {\rm
(\ref{eta-def})}.
  \item[{\rm (ii)}] {\rm [Extremals]} Let $x_0\in M$. The following
statements are equivalent:
\begin{itemize}
  \item[\rm (a)] For every $\kappa>0$ there exists a non-negative extremal function $u\in {\rm Lip}_0(M)$
  in ${\bf ({MS})}_{{\sf C}_2(p,n)}^2$, concentrated around $x_0$
  and $\mathcal H^n({\rm sprt}\ u )=\kappa;$
   \item[\rm (b)] $(M,g)$ is isometric to $(\mathbb R^n,e)$.
\end{itemize}
\end{itemize}
\end{theorem}

\begin{theorem}\label{theorem-morrey-4} Let $(M,g)$ be a  complete,
 $n-$di\-men\-sional Riemannian manifold with non-negative Ricci
 curvature, let $p>n,$ and assume that ${\bf ({MS})}_{C}^2$ holds on $(M,g)$
 for some $C>0.$ Then the following assertions hold:
 \begin{itemize}
   \item[\rm (i)] $C\geq {\sf C}_2(p,n);$
   \item[\rm (ii)] $(M,g)$ has the large volume balls property, i.e., $${\rm
   Vol}_g(B(x,\rho))\geq \left(\frac{{\sf C}_2(p,n)}{C}\right)^{\frac{pn}{p-n}+1}\omega_n\rho^n\ {for\ all}\ x\in M,\rho\geq 0;$$
   \item[\rm (iii)] ${\bf ({MS})}_{{\sf C}_2(p,n)}^2$ holds on
   $(M,g)$ if and only if $(M,g)$ is isometric to $(\mathbb R^n,e).$
 \end{itemize}
\end{theorem}


{\it Organization of the paper.} In Section \ref{sect-2} we recall
the notions and results from Riemannian geometry which are used
throughout the proofs. In Section \ref{sect-3} we  deal with the
sharp Morrey-Sobolev inequality  with support-bound, providing the
proof of Theorems \ref{theorem-morrey-1} and \ref{theorem-morrey-2}.
In Section \ref{sect-4} we treat the sharp Morrey-Sobolev inequality
with $L^1$-bound, proving Theorems \ref{theorem-morrey-3} and
\ref{theorem-morrey-4}.


\section{Preliminaries}\label{sect-2}

Let $(M,g)$ be an complete $n-$dimensional  Riemannian manifold, and
$d:M\times M\to [0,\infty)$ be the metric function associated to the
Riemannian metric $g$. Let $B(x,\rho)=\{y\in M:d(x,y)<\rho\}$ be the
open geodesic ball with center $x\in M$ and radius $\rho>0.$ If
 ${\text d}V_g$ is the canonical volume
element on $(M,g)$, the volume of an open bounded set $S\subset M$
is Vol$_g(S)=\int_S {\text d}V_g=\mathcal H^n(S)$, where $\mathcal
H^n(S)$ is the $n-$dimensional Hausdorff measure of $S$ with respect
to the metric function $d$. If ${\text d}\sigma_g$ denotes the
$(n-1)-$dimensional Riemann measure induced  on $\partial S$ by $g$,
 Area$_g(\partial S)=\int_{\partial S} {\text d}\sigma_g=\mathcal
H^{n-1}(\partial S)$ denotes the area of $\partial S$ with respect
to the metric $g$. In general, one has for every $x\in M$ that
\begin{equation}\label{volume-comp-nullaban}
   \lim_{\rho\to 0^+}\frac{{\rm Vol}_g(B(x,\rho))}{\omega_n
   \rho^n}=1,
\end{equation}
where $\omega_n$ is the volume of the standard $n-$dimensional
Euclidean unit ball. As usual, $B_e(0,\delta)$, ${\rm d}x$,  ${\rm
d}\sigma_e$, Vol$_e(S)$ and Area$_e(S)$ denote the Euclidean
counterparts of the above notions when  $S\subset \mathbb R^n$.

Let $p>1.$ The norm of $L^p(M)$ is given by
$\|u\|_{L^p(M)}=\left(\displaystyle\int_M |u|^p{\rm
d}V_g\right)^\frac{1}{p}$. Let $u:M\to \mathbb R$ be a function of
class $C^1.$ If $(x^i)$ denotes the local coordinate system on a
coordinate neighborhood of $x\in M$, and the local components of the
differential of $u$ are denoted by $u_i=\frac{\partial u}{\partial
x_i}$, then the local components of the gradient  $\nabla_g u$ are
$u^i=g^{ij}u_j$. Here, $g^{ij}$ are the local components of
$g^{-1}=(g_{ij})^{-1}$. The $L^p(M)$ norm of  $\nabla_g u(x)\in
T_xM$ is given by $\|\nabla_g u\|_{L^p(M)}=\left(\displaystyle\int_M
|\nabla_gu|^p{\rm d}V_g\right)^\frac{1}{p}.$ If  $u\in{\rm
Lip}_0(M)$, i.e., $u:M\to \mathbb R$ is a Lipschitz function with
compact support, then it is a.e. differentiable on $M$ and
$\|\nabla_g u\|_{L^p(M)}$ is well-defined. The space $W^{1,p}(M)$ is
the completion of $C_0^\infty(M)$ w.r.t. the norm
$\|u\|_{W^{1,p}(M)}=\|u\|_{L^p(M)}+\|\nabla_g u\|_{L^p(M)}.$

In the proof of our results Bishop-Gromov-type volume comparison
principles  play a crucial role. On account of Wu and Xin
\cite[Theorems 6.1 \& 6.3]{Wu-Xin}, we adapt the following version:

\begin{theorem}\label{comparison-volume}{\rm [Volume comparison]}
Let $(M,g)$ be a complete, $n-$di\-men\-sional Riemannian manifold
and $x_0\in M.$ Then the following statements hold.
\begin{itemize}
  \item[{\rm (a)}] If $(M,g)$ is a Cartan-Hadamard manifold,
   the function
$\rho\mapsto \frac{{\rm Vol}_g(B(x_0,\rho))}{\rho^n}$ is
non-decreasing, $\rho>0$. 
In particular, from {\rm (\ref{volume-comp-nullaban})} we have
\begin{equation}\label{volume-comp-altalanos-0}
{{\rm Vol}_g(B(x_0,\rho))}\geq \omega_n \rho^n\ {for\ all} \ \rho>0.
\end{equation}
If equality holds in {\rm (\ref{volume-comp-altalanos-0})}, then the
sectional curvature is identically zero.
   \item[{\rm (b)}] If $(M,g)$ has non-negative Ricci curvature,
   the function
$\rho\mapsto \frac{{\rm Vol}_g(B(x_0,\rho))}{\rho^n}$ is
non-increasing, $\rho>0$. In particular, from {\rm
(\ref{volume-comp-nullaban})} we have
\begin{equation}\label{volume-comp-altalanos-2}
{{\rm Vol}_g(B(x_0,\rho))}\leq \omega_n \rho^n\ {for\ all}\ \
\rho>0.
\end{equation}
If equality holds in {\rm (\ref{volume-comp-altalanos-2})}, then the
sectional curvature is identically zero.
\end{itemize}
\end{theorem}

\section{Sharp Morrey-Sobolev inequality  with
support-bound}\label{sect-3}

Let $(M,g)$ be a complete $n-$dimensional Riemannian manifold, and
let $p>n.$ For $C>0,$ we recall the Morrey-Sobolev inequality ${\bf
({MS})}_{C}^1$ with support-bound,
 i.e.,
 {$$
    \| u\|_{L^\infty(M)} \leq C \mathcal H^n({\rm sprt}\ u )^{\frac{1}{n}-\frac{1}{p}} \| \nabla_gu\|_{L^p(M)},\ \ \forall u\in {\rm Lip}_0(M).
$$}
We first present a result inspired by Aubin \cite{Aubin} and Hebey
\cite{Hebey}.

\begin{proposition}\label{prop-elso} If ${\bf ({MS})}_{C}^1$ holds, then $C\geq {\sf C}_1(p,n).$
\end{proposition}

{\it Proof.} Assume by contradiction that $C< {\sf C}_1(p,n).$ Let
$x_0\in M$.  For every $\varepsilon>0$, there exists a local chart
$(\Omega,\phi)$ of $M$ at the point $x_0$ and a number $\delta>0$
such that $\phi(\Omega)=B_e(0,\delta)$ and the components $g_{ij}$
of the metric $g$ satisfy
\begin{equation}\label{two-sided}
    (1-\varepsilon)\delta_{ij}\leq g_{ij} \leq (1+\varepsilon)\delta_{ij}
\end{equation}
in the sense of bilinear forms.

Due to ${\bf ({MS})}_{C}^1$  and to the two-sided estimate
(\ref{two-sided}), for $\varepsilon>0$ small enough, there exists
$\tilde \delta>0$ and $C'< {\sf C}_1(p,n)$ such that for every
$\delta\in (0,\tilde \delta)$ and $w\in {\rm Lip}_0(B_e(0,\delta))$,
\begin{equation}\label{c-hpw-uj-meg}
     \| w\|_{L^\infty(B_e(0,\delta))} \leq C' \mathcal H^n({\rm sprt}\ w )^{\frac{1}{n}-\frac{1}{p}} \| \nabla w\|_{L^p(B_e(0,\delta))}.
\end{equation}
Let $u\in {\rm Lip}_0(\mathbb R^n)$ be arbitrarily fixed and set
$w_\lambda(x)=u(\lambda x)$, $\lambda>0.$ For enough large
$\lambda>0,$ one has  $w_\lambda\in {\rm Lip}_0(B_e(0,\delta))$.
Replacing $w_\lambda$ into (\ref{c-hpw-uj-meg}), and using the
scaling properties
$$\| w_\lambda\|_{L^\infty(B_e(0,\delta))}=\| u\|_{L^\infty(\mathbb R^n)},\ \ \ \mathcal H^n({\rm sprt}\ w_\lambda )=\lambda^{-n}\mathcal H^n({\rm sprt}\ u ),$$
and
$$\int_{B_e(0,\delta)}|\nabla w_\lambda|^p{\text d}x=\lambda^{p-n}\int_{\mathbb R^n}|\nabla u|^p{\text d}x,$$
 one has
 $$ \| u\|_{L^\infty(\mathbb R^n)} \leq C' \mathcal H^n({\rm sprt}\ u )^{\frac{1}{n}-\frac{1}{p}} \| \nabla u\|_{L^p(\mathbb R^n)}.$$
 Inserting the function introduced in (\ref{extremal-1}) into the latter relation ,
 we obtain
   ${\sf C}_1(p,n)\leq
 C',$ a contradiction. \hfill $\square$\\

{\it Proof of Theorem \ref{theorem-morrey-1}.} The first part of the
proof is similar to Druet, Hebey and Vaugon \cite{Druetetal}, see
also Aubin, Druet and Hebey \cite{ADH}; since some intermediate
steps will be crucial in the second part (i.e., in the existence of
extremal functions), we shall present its complete proof. Let $p>n.$

(i)  Clearly, it is enough to consider only non-negative test
functions in the Morrey-Sobolev inequality ${\bf ({MS})}_{C}^1$.
Moreover, by standard approximation/density argument and Morse
theory, it is sufficient to deal with continuous
 test functions $u:M\to [0,\infty)$ having compact support $S\subset
 M$, where $S$ is an enough smooth set, $u$ being of class
 $C^2$ in $S$ and having only non-degenerate critical points in $S.$
 Fixing such a function $u:M\to [0,\infty)$, we associate to $u$ its
 Euclidean decreasing rearrangement function $u^*:\mathbb R^n\to
 [0,\infty)$ which is
 radially symmetric and is defined  for every $t>0$ by
 \begin{equation}\label{vol-egyenloseg}
    {\rm Vol}_e(\{x\in \mathbb R^n:u^*(x)>t\})={\rm Vol}_g(\{x\in
    M:u(x)>t\})\stackrel{\rm def.}{=:}V(t).
 \end{equation}
 By definition, $u^*$ is a Lipschitz
function with compact support, and
\begin{equation}\label{sup-ket-fuggveny}
    \|u\|_{L^\infty(M)}=\|u^*\|_{L^\infty(\mathbb R^n)},\ \ \ {\rm Vol}_g({\rm sprt}\ u
    )={\rm Vol}_e({\rm sprt}\ u^* ).
\end{equation}

On one hand, for every $0<t<\|u\|_{L^\infty(M)},$ we consider the
level sets
$$\Gamma_t=u^{-1}(t)\subset S\subset M,\ \ \ \ \Gamma_t^*=(u^*)^{-1}(t)\subset \mathbb R^n.$$
Since $u^*$ is radial, $\Gamma_t^*$ is an $(n-1)-$dimensional sphere
with $\Gamma_t^*=\partial (\{x\in \mathbb R^n:u^*(x)>t\})$ for every
$0<t<\|u\|_{L^\infty(M)},$ and
$${\rm Area}_e(\Gamma_t^*) =  n\omega_n^\frac{1}{n} {\rm
Vol}_e(\{x\in \mathbb R^n:u^*(x)>t\})^\frac{n-1}{n}.$$ In
particular,  the latter relation, the validity of Cartan-Hadamard
conjecture and (\ref{vol-egyenloseg}) imply that
\begin{equation}\label{areak-egyenlotlenseg}
    {\rm Area}_g(\Gamma_t) \geq  {\rm Area}_e(\Gamma_t^*)\ \ {\rm for\ every}\ 0<t<\|u\|_{L^\infty(M)}.
\end{equation}
A simple application of the co-area formula (see Chavel \cite[p.
302]{Chavel}) and (\ref{vol-egyenloseg}) give
\begin{equation}\label{V-deriv}
    V'(t)=-\int_{\Gamma_t}\frac{1}{|\nabla_g u|}{\rm d
\sigma}_g=-\int_{\Gamma_t^*}\frac{1}{|\nabla u^*|}{\rm d \sigma}_e.
\end{equation}
Since $|\nabla u^*|$ is constant on the sphere $\Gamma_t^*$, the
second relation from (\ref{V-deriv}) gives that
\begin{equation}\label{V-DER-MASIK}
V'(t)=-\frac{{\rm Area}_e(\Gamma_t^*)}{|\nabla u^*(x)|},\ x\in
\Gamma_t^*.
\end{equation}

On the other hand, by H\"older's inequality and the first relation
of (\ref{V-deriv}), one has
$${\rm Area}_g(\Gamma_t)=\int_{\Gamma_t}{\rm d \sigma}_g\leq\left(-V'(t)\right)^\frac{p-1}{p}\left(\int_{\Gamma_t}{|\nabla_g u|^{p-1}}{\rm d
\sigma}_g\right)^\frac{1}{p}.$$ Consequently, by
(\ref{areak-egyenlotlenseg}) and (\ref{V-DER-MASIK}), for every
$0<t<\|u\|_{L^\infty(M)}$ we have
\begin{eqnarray*}
  \int_{\Gamma_t}{|\nabla_g u|^{p-1}}{\rm d
\sigma}_g &\geq& {\rm Area}_g(\Gamma_t)^p \left(-V'(t)\right)^{1-p}\\
   &\geq & {\rm Area}_e(\Gamma_t^*)^p \left(\frac{{\rm Area}_e(\Gamma_t^*)}{|\nabla u^*(x)|}\right)^{1-p}\ \ \ \ \ \ \ \ (x\in \Gamma_t^*)\\
   &=& \int_{\Gamma_t^*}{|\nabla u^*|^{p-1}}{\rm d
\sigma}_e.
\end{eqnarray*}
By the co-area formula and the latter inequality, an integration
with respect to $t$ gives
\begin{equation}\label{Polya-Szego}
    \int_{M}{|\nabla_g u|^{p}}{\rm d
}V_g\geq \int_{\mathbb R^n}{|\nabla u^*|^{p}}{\rm d}x.
\end{equation}
Applying  Talenti's inequality for the function $u^*:\mathbb R^n\to
\mathbb R$ (see \cite[Theorem 2.E]{Talenti} and Introduction),
relations (\ref{sup-ket-fuggveny}) and (\ref{Polya-Szego}) provide
\begin{eqnarray}\label{utolso-becs}
 \nonumber \|u\|_{L^\infty(M)} &=& \|u^*\|_{L^\infty(\mathbb R^n)} \\
   &\leq & {\sf C}_1(p,n) \mathcal H^n({\rm sprt}\ u^* )^{\frac{1}{n}-\frac{1}{p}} \| \nabla u^*\|_{L^p(\mathbb R^n)} \nonumber \\
   &\leq & {\sf C}_1(p,n) \mathcal H^n({\rm sprt}\ u )^{\frac{1}{n}-\frac{1}{p}} \|
   \nabla_g
   u\|_{L^p(M)},
\end{eqnarray}
which is precisely ${\bf ({MS})}_{{\sf C}_1(p,n)}^1$  on $(M,g).$
Moreover, this inequality and Proposition \ref{prop-elso} show that
$${\sf C}_1(p,n)^{-1}=\inf_{u\in {\rm Lip}_0(M)\setminus \{0\}} \frac{\mathcal H^n({\rm sprt}\ u )^{\frac{1}{n}-\frac{1}{p}} \| \nabla_gu\|_{L^p(M)}}{\|
u\|_{L^\infty(M)}}.$$

(ii)  By Talenti's result, we clearly have (b)$\Rightarrow$(a). Let
$x_0\in M,$ and assume that (a) holds, i.e., for every $\kappa>0$
there exists a non-negative extremal function $u\in {\rm Lip}_0(M)$
  in ${\bf ({MS})}_{{\sf C}_1(p,n)}^1$, concentrated around $x_0$
  and $\mathcal H^n({\rm sprt}\ u )=\kappa.$ Therefore, in
  (\ref{utolso-becs}) we have equalities; in particular, the Euclidean decreasing rearrangement function $u^*:\mathbb R^n\to
 [0,\infty)$ associated to $u$ verifies $$\|u^*\|_{L^\infty(\mathbb R^n)}={\sf C}_1(p,n) \mathcal H^n({\rm sprt}\ u^* )^{\frac{1}{n}-\frac{1}{p}} \| \nabla u^*\|_{L^p(\mathbb R^n)}.$$
According to Talenti's result, since $u^*$ is an extremal in ${\bf
({MS})}_{{\sf C}_1(p,n)}^1$ on $(\mathbb R^n,e)$,  it has the shape
from (\ref{extremal-1}), i.e.,
$$ u^*(x)=\left(\lambda^\frac{p-n}{p-1}-|x|^\frac{p-n}{p-1}\right)_+,\
x\in \mathbb R^n.$$ In addition, since $\mathcal H^n({\rm sprt}\ u^*
)=\mathcal H^n({\rm sprt}\ u)=\kappa,$  we have
$\lambda=(\kappa\omega_n^{-1})^\frac{1}{n}.$ Consequently, for every
$0<t<\|u\|_{L^\infty(M)}=\|u^*\|_{L^\infty(\mathbb
R^n)}=\lambda^\frac{p-n}{p-1}$, one has $\{x\in \mathbb
R^n:u^*(x)>t\}=B_e(0,\rho_t),$ where
$$\rho_t=\left(\lambda^\frac{p-n}{p-1}-t\right)^\frac{p-1}{p-n}.$$

Fix $0<t<\lambda^\frac{p-n}{p-1}.$ We claim that
$$\{x\in M:u(x)>t\}=B(x_0,\rho_t).$$
On one hand, since $u$ is concentrated around $x_0$, there exists
$\rho_t'>0$ such that $\{x\in M:u(x)>t\}=B(x_0,\rho_t').$ Therefore,
the claim is concluded once we prove that $\rho_t'=\rho_t.$ Due to
(\ref{vol-egyenloseg}), one has
\begin{equation}\label{iso-utolsok}
    {\rm Vol}_g(B(x_0,\rho_t'))={\rm Vol}_e(B_e(0,\rho_t)).
\end{equation}
On  the other hand, since $u$ is an extremal in ${\bf ({MS})}_{{\sf
C}_1(p,n)}^1$, we have equalities not only in (\ref{utolso-becs})
but also in (\ref{Polya-Szego}). Subsequently, we have equality also
in (\ref{areak-egyenlotlenseg}), i.e.,
    $${\rm Area}_g(\Gamma_t) =  {\rm Area}_e(\Gamma_t^*).$$
This relation together with (\ref{vol-egyenloseg}) imply that we
have equality case in the Cartan-Hadamard conjecture, i.e., $\{x\in
M:u(x)>t\}=B(x_0,\rho_t')$ is isometric to the $n-$dimensional
Euclidean ball with volume ${\rm Vol}_g(B(x_0,\rho_t'))$. On account
of (\ref{iso-utolsok}), we actually have that $B(x_0,\rho_t')$ and
$B_e(0,\rho_t)$ are isometric, thus $\rho_t'=\rho_t.$ Therefore,
${\rm Vol}_g(B(x_0,\rho_t))=\omega_n\rho_t^n.$ If $t\to 0^+$, the
latter relation implies that ${\rm
Vol}_g(B(x_0,\lambda))=\omega_n\lambda^n.$ Due to the arbitrariness
of $\kappa>0$, so $\lambda=(\kappa\omega_n^{-1})^\frac{1}{n},$ we
have that
$${\rm Vol}_g(B(x_0,\rho))=\omega_n\rho^n\ \ {\rm for\ all}\ \rho>0.$$
By Theorem \ref{comparison-volume} (i) we have that the sectional
curvature on the Cartan-Hadamard manifold $(M,g)$ is identically
zero, which concludes the proof.\hfill $\square$\\

{\it Proof of Theorem \ref{theorem-morrey-2}.} We assume that ${\bf
({MS})}_{C}^1$ holds on $(M,g)$ for some $C>0.$ By Proposition
\ref{prop-elso}, we already have that $C\geq {\sf C}_1(p,n),$ i.e.,
(i) is proved.

(ii) We are going to prove that  $(M,g)$ has the large volume balls
property, i.e., $${\rm
   Vol}_g(B(x,\rho))\geq \left(\frac{{\sf C}_1(p,n)}{C}\right)^\frac{pn}{p-n}\omega_n\rho^n\ {\rm for\ all}\ x\in M,\rho\geq 0.$$
Let $x_0\in M$ be fixed. For every $\lambda>0,$  we consider the
function
$$
u_\lambda(x)=\left(\lambda^\frac{p-n}{p-1}-d(x_0,x)^\frac{p-n}{p-1}\right)_+,\
x\in M.$$ It is clear that $u_\lambda\in {\rm Lip}_0(M)$ and
\begin{equation}\label{tulajd-morrey-1}
    \|u_\lambda\|_{L^\infty(M)}=\lambda^\frac{p-n}{p-1},\ \ \ \mathcal H^n({\rm sprt}\
    u_\lambda
    )={\rm Vol}_g(B(x_0,\lambda)).
\end{equation}
The chain rule (see Hebey \cite[Proposition 2.5]{Hebey}) implies
that $$\nabla_g
u_\lambda(x)=-\frac{p-n}{p-1}d(x_0,x)^\frac{1-n}{p-1}\nabla_g
d(x_0,x),\ \ x\in B(x_0,\lambda).$$ Taking into account that
$|\nabla_g d(x_0,x)|=1$ for a.e. $x\in M$, the layer cake
representation and Theorem \ref{comparison-volume} (ii) give that
\begin{eqnarray*}
  \int_M |\nabla_g
u_\lambda|^p{\rm d}V_g &=& \int_{B(x_0,\lambda)} |\nabla_g
u_\lambda|^p{\rm d}V_g \\
   &=& \left(\frac{p-n}{p-1}\right)^p\int_{B(x_0,\lambda)}d(x_0,x)^\frac{p(1-n)}{p-1}{\rm d}V_g \\
   &=& \left(\frac{p-n}{p-1}\right)^p\int_0^\infty {\rm Vol}_g\left(\{x\in B(x_0,\lambda):d(x_0,x)^\frac{p(1-n)}{p-1}>t\}\right){\rm
   d}t\\
 &=& \left(\frac{p-n}{p-1}\right)^p\int_{\lambda^\frac{p(1-n)}{p-1}}^\infty {\rm Vol}_g\left(\{x\in B(x_0,\lambda):d(x_0,x)^\frac{p(1-n)}{p-1}>t\}\right){\rm d}t \\
    && +\left(\frac{p-n}{p-1}\right)^p\int_0^{\lambda^\frac{p(1-n)}{p-1}}{\rm Vol}_g\left(\{x\in B(x_0,\lambda):d(x_0,x)^\frac{p(1-n)}{p-1}>t\}\right){\rm d}t \\
&=& \left(\frac{p-n}{p-1}\right)^p\frac{p(n-1)}{p-1}\int_0^\lambda {\rm Vol}_g\left( B(x_0,\rho)\right)\rho^\frac{-pn+1}{p-1}{\rm d}\rho \\
&& +\left(\frac{p-n}{p-1}\right)^p {\rm Vol}_g\left( B(x_0,\lambda)\right)\lambda^\frac{p(1-n)}{p-1}\\
&\leq& \left(\frac{p-n}{p-1}\right)^p\omega_n\left[\frac{p(n-1)}{p-1}\int_0^\lambda \rho^{n+\frac{-pn+1}{p-1}}{\rm d}\rho +\lambda^{n+\frac{p(1-n)}{p-1}}\right]\\
&=&\left(\frac{p-n}{p-1}\right)^{p-1}n\omega_n
\lambda^\frac{p-n}{p-1}.
\end{eqnarray*}

Inserting $u_\lambda$ into  ${\bf ({MS})}_{C}^1$, relation
(\ref{tulajd-morrey-1}) and the above estimate yield that
$$\lambda^\frac{p-n}{p-1}\leq C{\rm
Vol}_g(B(x_0,\lambda))^{\frac{1}{n}-\frac{1}{p}}\left(\frac{p-n}{p-1}\right)^\frac{1}{p'}(n\omega_n)^\frac{1}{p}
\lambda^\frac{p-n}{p(p-1)}.$$ Reorganizing this inequality and
taking into account the form of the constant  ${\sf C}_1(p,n)$, see
(\ref{C1-konstans}), it turns out that for every $\lambda>0$, we
have
\begin{equation}\label{x0-as-ossz}
{\rm
   Vol}_g(B(x_0,\lambda))\geq \left(\frac{{\sf C}_1(p,n)}{C}\right)^\frac{pn}{p-n}\omega_n\lambda^n.
\end{equation}

Let $x\in M$ and $\rho>0$ be fixed arbitrarily. Then,  one has
\begin{eqnarray*}
   \frac{{\rm
   Vol}_g(B(x,\rho))}{\omega_n\rho^n} &\geq& \limsup_{r\to \infty}\frac{{\rm
   Vol}_g(B(x,r))}{\omega_nr^n}\ \ \ \ \ \ \ \ \ \ \ \ \  \ \ \ \ \ \ \ \   {\rm [cf.\
   Theorem\
\ref{comparison-volume}\ (ii)]} \\
   &\geq&   \limsup_{r\to \infty}\frac{{\rm
   Vol}_g(
    B(x_0,r-{d}(x_0,x)))}{\omega_nr^n}\ \ \ \ \ \  {\rm [}B(x,r)\supset B(x_0,r-{d}(x_0,x)){\rm ]} \\
   &=& \limsup_{r\to \infty}\left(\frac{{\rm
   Vol}_g(
    B(x_0,r-{d}(x_0,x)))}{\omega_n(r-{d}(x_0,x))^n}\cdot\frac{(r-{d}(x_0,x))^n}{r^n}\right)
\\&\geq&\left(\frac{{\sf C}_1(p,n)}{C}\right)^\frac{pn}{p-n},\ \ \ \ \ \ \ \ \ \ \
\ \ \ \ \ \ \ \ \ \  \  \ \ \ \ \ {\rm [cf.\
 (\ref{x0-as-ossz})]}
\end{eqnarray*}
which  concludes the proof of (ii).

(iii) By Talenti's result, if $(M,g)$ is isometric to $(\mathbb
R^n,e),$ then ${\bf ({MS})}_{{\sf C}_1(p,n)}^1$ holds. Conversely,
let us assume that ${\bf ({MS})}_{{\sf C}_1(p,n)}^1$ holds on
$(M,g)$. First, by (ii) we have that ${\rm
   Vol}_g(B(x_,\rho))\geq \omega_n\rho^n$ for every $x\in M$ and $\rho>0.$
By (\ref{volume-comp-altalanos-2}), we also have  the converse
inequality ${\rm
   Vol}_g(B(x,\rho))\leq \omega_n\rho^n,$ thus
\begin{equation}\label{lllll}
    {\rm
   Vol}_g(B(x_,\rho))= \omega_n\rho^n\ \ {\rm for\ all}\ x\in M,\ \rho>0.
\end{equation}
By Theorem \ref{comparison-volume} (ii), it follows that the
sectional curvature on $(M,g)$ is identically zero. Then relation
(\ref{lllll})  implies that $(M,g)$ is isometric to $(\mathbb
R^n,e).$\hfill $\square$\\

\section{Sharp
Morrey-Sobolev inequality with $L^1$-bound}\label{sect-4}

The structure of this section is similar to the previous one; in the
sequel, we shall point out the differences. As before, let $(M,g)$
be a complete $n-$dimensional Riemannian manifold, and let $p>n.$
For $C>0,$ we recall the Morrey-Sobolev inequality  with
$L^1-$bound,
 i.e.,
 $$
    \| u\|_{L^\infty(M)} \leq C \|u\|_{L^1(M)}^{1-\eta} \| \nabla_gu\|_{L^p(M)}^\eta,\ \ \forall u\in {\rm Lip}_0(M), \eqno{{\bf ({MS})}_{C}^2}
$$
where
$$
    \eta=\frac{np}{np+p-n}.
$$

\begin{proposition}\label{prop-masodik} If ${\bf ({MS})}_{C}^2$ holds, then $C\geq {\sf C}_2(p,n).$
\end{proposition}

{\it Proof.}  We follow the proof of Proposition \ref{prop-elso}.
The only minor difference is that we use a further scaling property.
Namely, let
 $u\in {\rm Lip}_0(\mathbb R^n)$  and  $w_\lambda(x)=u(\lambda x)$,
 $\lambda>0.$ Then
for enough large $\lambda>0,$ one has  $w_\lambda\in {\rm
Lip}_0(B_e(0,\delta))$, and in addition to the scaling properties
from Proposition \ref{prop-elso} we also have $ \| w_\lambda
\|_{L^1(\mathbb R^n)}=\lambda^{-n}\| u \|_{L^1(\mathbb R^n)}.$
 \hfill $\square$\\

{\it Proof of Theorem \ref{theorem-morrey-3}.} Let $p>n.$ Let
$u:M\to [0,\infty)$ be a function with the same properties as in the
proof of Theorem \ref{theorem-morrey-1}. If we associate to $u$ its
 Euclidean decreasing rearrangement function $u^*:\mathbb R^n\to
 [0,\infty)$ which is
 radially symmetric and defined  by
(\ref{vol-egyenloseg}), one has that
$$\|u\|_{L^\infty(M)}=\|u^*\|_{L^\infty(\mathbb R^n)},\ \  \|\nabla_g u\|_{L^p(M)}\geq \|\nabla u^*\|_{L^p(\mathbb R^n)},$$
see (\ref{sup-ket-fuggveny}) and (\ref{Polya-Szego}), respectively.
In addition, by the layer cake representation and
(\ref{vol-egyenloseg}), we also have
\begin{eqnarray*}
  \|u\|_{L^1(M)} &=& \int_0^\infty {\rm Vol}_g(\{x\in M:u(x)>t\}){\rm
d}t=\int_0^\infty{\rm Vol}_e(\{x\in \mathbb R^n:u^*(x)>t\}){\rm d}t \\
   &=& \|u^*\|_{L^1(\mathbb R^n)}.
\end{eqnarray*}
Consequently, the latter relations and Talenti's result (see
\cite[Theorem 2.C]{Talenti} and the Introduction) imply that
\begin{eqnarray}\label{utolso-becs-masik}
 \nonumber \|u\|_{L^\infty(M)} &=& \|u^*\|_{L^\infty(\mathbb R^n)} \\
   &\leq & {\sf C}_2(p,n) \|u^*\|_{L^1(\mathbb R^n)}^{1-\eta} \| \nabla u^*\|_{L^p(\mathbb R^n)}^\eta \nonumber \\
   &\leq & {\sf C}_2(p,n) \|u\|_{L^1(M)}^{1-\eta} \|
   \nabla_gu\|_{L^p(M)}^\eta,
\end{eqnarray}
i.e., ${\bf ({MS})}_{{\sf C}_2(p,n)}^2$ holds on $(M,g).$ The
sharpness of the constant ${\sf C}_2(p,n)$ follows by the latter
estimate and Proposition \ref{prop-masodik}, concluding the proof of
(i).

Before to provide the proof of (ii), we introduce some notations
which will be useful in the sequel. Let $\lambda>0$ and define
$f_\lambda:(0,\lambda]\to [0,\infty)$ and $F_\lambda:[0,\lambda]\to
[0,\infty)$ by
\begin{equation}\label{f-lambda}
f_\lambda(r)=r^\frac{1-n}{p-1}(\lambda^n-r^n)^\frac{1}{p-1},\ \
F_\lambda(s)=\int_0^s f_\lambda(r) {\rm d}r.
\end{equation}

(ii) Let $x_0\in M$ and $\kappa>0$ be fixed arbitrarily and let
$u\in {\rm Lip}_0(M)$ be a non-negative extremal function in ${\bf
({MS})}_{{\sf C}_2(p,n)}^2$,  concentrated around $x_0$
  and $\mathcal H^n({\rm sprt}\ u )=\kappa.$ Since we have equalities in
  (\ref{utolso-becs-masik}), the function $u^*:\mathbb R^n\to
 [0,\infty)$  verifies $$\|u^*\|_{L^\infty(\mathbb R^n)}={\sf C}_2(p,n) \|u^*\|_{L^1(\mathbb R^n)}^{1-\eta} \| \nabla u^*\|_{L^p(\mathbb R^n)}^\eta.$$
Since $u^*$ is an extremal in ${\bf ({MS})}_{{\sf C}_2(p,n)}^2$ on
$(\mathbb R^n,e)$, by Talenti's result, its expression is given by
$$u^*(x)=\left\{
\begin{array}{lll}
F_\lambda(\lambda)-F_\lambda(|x|), & & {\rm if}\ x\in
B_e(0,\lambda);
\\ 0, & & {\rm if}\ x\notin
B_e(0,\lambda).
\end{array}
\right.$$ Note that by $\mathcal H^n({\rm sprt}\ u^* )=\mathcal
H^n({\rm sprt}\ u)=\kappa,$  we have
$\lambda=(\kappa\omega_n^{-1})^\frac{1}{n}.$ Moreover,
$$\|u\|_{L^\infty(M)}=\|u^*\|_{L^\infty(\mathbb
R^n)}=F_\lambda(\lambda)==\frac{\lambda^{p'}}{n}{\sf
    B}\left(\frac{1-n}{n}p'+1,p'\right)$$ and since $F_\lambda$ is increasing on
$[0,\lambda],$ for every $0<t<F_\lambda(\lambda)$ one has $\{x\in
\mathbb R^n:u^*(x)>t\}=B_e(0,\rho_t),$ where
$$\rho_t=F_\lambda^{-1}(F_\lambda(\lambda)-t).$$
A similar argument as in the proof of Theorem \ref{theorem-morrey-1}
(ii) shows that for every  $0<t<F_\lambda(\lambda),$ we have $\{x\in
M:u(x)>t\}=B(x_0,\rho_t),$ and finally
$${\rm Vol}_g(B(x_0,\rho))=\omega_n\rho^n\ \ {\rm for\ all}\ \rho>0,$$
which concludes the proof.  \hfill $\square$\\

{\it Proof of Theorem \ref{theorem-morrey-4}.} (i) By Proposition
\ref{prop-masodik}, we have that $C\geq {\sf C}_2(p,n)$ whenever
${\bf ({MS})}_{C}^2$ is assumed to hold on $(M,g)$.

(ii) Let $x_0\in M$ be fixed. By using (\ref{f-lambda}), for every
$\lambda>0$ we consider the functions $u_\lambda\in {\rm Lip}_0(M)$
and $w_\lambda\in {\rm Lip}_0(\mathbb R^n)$ defined by
$$
u_\lambda(x)=\left\{
\begin{array}{lll}
F_\lambda(\lambda)-F_\lambda(d(x_0,x)), & & {\rm if}\ x\in
B(x_0,\lambda);
\\ 0, & & {\rm if}\ x\notin
B(x_0,\lambda),
\end{array}
\right.$$ and
$$
w_\lambda(x)=\left\{
\begin{array}{lll}
F_\lambda(\lambda)-F_\lambda(|x|), & & {\rm if}\ x\in
B_e(0,\lambda);
\\ 0, & & {\rm if}\ x\notin
B_e(0,\lambda).
\end{array}
\right.$$ Since $u_\lambda$ verifies  ${\bf ({MS})}_{C}^2$ on
$(M,g),$ and $w_\lambda$ is an extremal in  ${\bf ({MS})}_{{\sf
C}_2(p,n)}^2$ on $(\mathbb R^n,e)$, we have that
\begin{equation}\label{u-lambda-egy}
    \|u_\lambda\|_{L^\infty(M)}\leq C
\|u_\lambda\|_{L^1(M)}^{1-\eta} \| \nabla_g
u_\lambda\|_{L^p(M)}^\eta
\end{equation}
and
\begin{equation}\label{w-lambda-ident}
    \|w_\lambda\|_{L^\infty(\mathbb R^n)}={\sf C}_2(p,n)
\|w_\lambda\|_{L^1(\mathbb R^n)}^{1-\eta} \| \nabla
w_\lambda\|_{L^p(\mathbb R^n)}^\eta.
\end{equation}
Moreover, by the above definitions and a simple computation give
\begin{equation}\label{tulajd-morrey-mmm}
    \|u_\lambda\|_{L^\infty(M)}=\|w_\lambda\|_{L^\infty(\mathbb
    R^n)}=F_\lambda(\lambda)=\frac{\lambda^{p'}}{n}{\sf
    B}\left(\frac{1-n}{n}p'+1,p'\right).
\end{equation}
 Similar
computations show that
\begin{equation}\label{w-lambda-1}
    \|w_\lambda\|_{L^1(\mathbb R^n)}=\omega_n\int_0^\lambda
\rho^nf_\lambda(\rho){\rm
   d}\rho=\frac{\lambda^{n+p'}\omega_n}{n}{\sf
    B}\left(\frac{1-n}{n}p'+2,p'\right)
\end{equation}
 and
$$\|\nabla w_\lambda\|_{L^p(\mathbb R^n)}={\lambda^\frac{n+p'}{p}\omega_n}^\frac{1}{p}{\sf
    B}\left(\frac{1-n}{n}p'+1,p'+1\right)^\frac{1}{p}.$$

In the sequel, we shall estimate $\|u_\lambda\|_{L^1(M)}$ and
$\|\nabla_gu_\lambda\|_{L^p(M)}$. First, by the layer cake
representation, one has
\begin{eqnarray*}
  \|u_\lambda\|_{L^1(M)} &=& \int_M u_\lambda(x){\rm d}V_g =\int_{B(x_0,\lambda)} (F_\lambda(\lambda)-F_\lambda(d(x_0,x))){\rm d}V_g\\
   &=& \int_0^\infty {\rm Vol}_g(\{x\in B(x_0,\lambda):F_\lambda(\lambda)-F_\lambda(d(x_0,x))>t\}){\rm d}t\\
   &&\ \ \ \ \ \ \ \ \ \ \ \ \  \ \ \ \ \ \ \ \ \ \ \ \ \ \ \ \ \ \ \ \  \ \ \ \ [{\rm change\ of\ var.}\ t=F_\lambda(\lambda)-F_\lambda(\rho)]\\
   &=& \int_0^\lambda {\rm Vol}_g(B(x_0,\rho))f_\lambda(\rho){\rm
   d}\rho.
\end{eqnarray*}
Then, since $\nabla_g u_\lambda(x)=-f_\lambda(d(x_0,x))\nabla_g
d(x_0,x)$ for every $ x\in B(x_0,\lambda),$ by Bishop-Gromov
comparison theorem (see relation (\ref{volume-comp-altalanos-2})),
one has
\begin{eqnarray*}
  \|\nabla_gu_\lambda\|_{L^p(M)}^p &=& \int_{B(x_0,\lambda)} f_\lambda(d(x_0,x))^p{\rm d}V_g\\
  &=&  \int_{B(x_0,\lambda)}\left(\lambda^n d(x_0,x)^{1-n}-d(x_0,x)\right)^{p'} {\rm d}V_g \\
   &=& \int_0^\infty {\rm Vol}_g(\{x\in B(x_0,\lambda):\left(\lambda^n d(x_0,x)^{1-n}-d(x_0,x)\right)^{p'}>t\}){\rm d}t\\
   &&\ \ \ \ \ \ \ \ \ \ \ \ \  \ \ \ \ \ \ \ \ \ \ \ \ \ \ \ \ \ \ \ \ \ \  \ \ \ \ \left[{\rm change\ of\ var.}\ t=\left(\lambda^n \rho^{1-n}-\rho\right)^{p'}\right]\\
   &=& p'\int_0^\lambda {\rm Vol}_g(B(x_0,\rho))\left(\lambda^n \rho^{1-n}-\rho\right)^{p'-1}\left((n-1)\lambda^n \rho^{-n}+1\right){\rm
   d}\rho\\
   &\leq&  p'\int_0^\lambda \left(\lambda^n \rho^{1-n}-\rho\right)^{p'-1}\left((n-1)\lambda^n +\rho^{n}\right){\rm
   d}\rho\\
   &=& \lambda^{n+p'}\omega_n{\sf
    B}\left(\frac{1-n}{n}p'+1,p'+1\right)\\
    &=& \|\nabla w_\lambda\|_{L^p(\mathbb R^n)}^p.
\end{eqnarray*}
Subtracting (\ref{w-lambda-ident}) from (\ref{u-lambda-egy}),
relations (\ref{tulajd-morrey-mmm}), (\ref{w-lambda-1}) and the
above computations give that for every $\lambda>0,$
\begin{equation}\label{eleg-jo}
    \int_0^\lambda \left({\rm Vol}_g(B(x_0,\rho))-\left(\frac{{\sf
C}_2(p,n)}{C}\right)^\frac{1}{1-\eta}\omega_n\rho^n\right)f_\lambda(\rho){\rm
   d}\rho\geq 0.
\end{equation}

 We claim
that
\begin{equation}\label{asymptotic}
    \ell_\infty^{x_0}:=\limsup_{\rho\to \infty}\frac{{\rm Vol}_g(B(x_0,\rho))}{\omega_n\rho^n}\geq \left(\frac{{\sf
C}_2(p,n)}{C}\right)^\frac{1}{1-\eta}.
\end{equation}
Assuming the contrary, there exists $\varepsilon_0>0$ such that for
some $\rho_0>0$,
$$\frac{{\rm Vol}_g(B(x_0,\rho))}{\omega_n\rho^n}\leq \left(\frac{{\sf
C}_2(p,n)}{C}\right)^\frac{1}{1-\eta}-\varepsilon_0,\ \forall
\rho\geq \rho_0.$$  By the latter inequality and (\ref{eleg-jo}),
for every $\lambda>\rho_0$ we obtain that
\begin{eqnarray*}
  0 &\leq&  \int_0^\lambda \left({\rm Vol}_g(B(x_0,\rho))-\left(\frac{{\sf
C}_2(p,n)}{C}\right)^\frac{1}{1-\eta}\omega_n\rho^n\right)f_\lambda(\rho){\rm
   d}\rho \\
   &\leq & \int_0^{\rho_0} {\rm Vol}_g(B(x_0,\rho))f_\lambda(\rho){\rm
   d}\rho-\varepsilon_0\omega_n\int_{\rho_0}^\lambda \rho^nf_\lambda(\rho){\rm
   d}\rho\\&&-\left(\frac{{\sf
C}_2(p,n)}{C}\right)^\frac{1}{1-\eta}\omega_n\int_0^{\rho_0}\rho^nf_\lambda(\rho){\rm
   d}\rho.
\end{eqnarray*}
Rearranging the above inequality, by (\ref{volume-comp-altalanos-2})
it follows that
\begin{equation}\label{nem-tudom-dejsze}
    \varepsilon_0\int_0^\lambda
   \rho^nf_\lambda(\rho){\rm
   d}\rho\leq\left(1-\left(\frac{{\sf
C}_2(p,n)}{C}\right)^\frac{1}{1-\eta}+\varepsilon_0\right)\int_0^{\rho_0}
   \rho^nf_\lambda(\rho){\rm
   d}\rho.
\end{equation}
  According to  (\ref{w-lambda-1}) and to the fact that $$\int_0^{\rho_0}
   \rho^nf_\lambda(\rho){\rm
   d}\rho\leq \rho_0^{n+\frac{1-n}{p-1}}\lambda^\frac{n}{p-1},$$ inequality (\ref{nem-tudom-dejsze}) implies
   $$\varepsilon_0\frac{\lambda^{n+p'}}{n}{\sf
    B}\left(\frac{1-n}{n}p'+2,p'\right)\leq\left(1-\left(\frac{{\sf
C}_2(p,n)}{C}\right)^\frac{1}{1-\eta}+\varepsilon_0\right)\rho_0^{n+\frac{1-n}{p-1}}\lambda^\frac{n}{p-1},\
\lambda>\rho_0.$$ Since $n+p'>\frac{n}{p-1},$ letting $\lambda\to
+\infty$ in the latter inequality, we arrive to a contradiction. The
proof of (\ref{asymptotic})  is complete.

By Theorem \ref{comparison-volume} (ii) and (\ref{asymptotic}), it
follows that for every $\rho>0$,
$$
\frac{{\rm Vol}_g(B(x_0,\rho))}{\omega_n\rho^n}\geq
\ell_\infty^{x_0}\geq\left(\frac{{\sf
C}_2(p,n)}{C}\right)^\frac{1}{1-\eta}.
$$
If $x\in M$ and $\rho>0$ are arbitrarily fixed, a similar argument
applies as in the last step of the proof of Theorem
\ref{theorem-morrey-2}, where the latter inequality takes the role
of  (\ref{x0-as-ossz}).

(iii) Similar to the proof of  Theorem
\ref{theorem-morrey-2} (iii).  \hfill $\square$\\

\end{document}